\documentclass[a4paper,12pt]{article}
\usepackage{amsmath,amssymb,amsthm,amsfonts}
\begin{document}

\newtheorem{theorem}{Theorem}[section]
 \newtheorem{corollary}[theorem]{Corollary}
 \newtheorem{lemma}[theorem]{Lemma}

\title{\Large\bf Rainbow number of matchings\\ in regular bipartite graphs
\footnote {Supported by NSFC, PCSIRT and the ``973" program.} }

\author{\small Xueliang Li$^1$ and Zhixia Xu$^{1,2}$\\
[2mm]
\small $^1$Center for Combinatorics and LPMC-TJKLC\\
\small Nankai University, Tianjin 300071, China. Email: lxl@nankai.edu.cn\\
\small $^2$College of Mathematics and System Sciences, Xinjiang University\\
\small Urumuqi, 830046, China. Email: irisxuzx@gmail.com}
\date{}
\maketitle

\begin{abstract}
Given a graph $G$ and a subgraph $H$ of $G$, let $rb(G,H)$ be the
minimum number $r$ for which any edge-coloring of $G$ with $r$
colors has a rainbow subgraph $H$. The number $rb(G,H)$ is called
the rainbow number of $H$ with respect to $G$. Denote $mK_2$ a
matching of size $m$ and $B_{n,k}$ a $k$-regular bipartite graph
with bipartition $(X,Y)$ such that $|X|=|Y|=n$ and $k\leq n$. In
this paper we give an upper and lower bound for
$rb(B_{n,k},mK_2)$, and show that for given $k$ and $m$, if $n$ is
large enough, $rb(B_{n,k},mK_2)$ can reach the lower bound. We
also determine the rainbow number of matchings in paths and
cycles.\\[2mm]
\end{abstract}

\vspace{10pt}

\section{\large Introduction}

We use Bondy and Murty \cite{bondy} for terminology and notations
not defined here and consider simple, finite graphs only.

The Ramsey problem asks for the optimal total number of colors
used on the edges of a graph without creating a monochromatic
subgraph. In anti-Ramsey problems, we are interested in
heterochromatic or rainbow subgraphs instead of monochromatic
subgraphs in edge-colorings. Given a graph $G$ and a subgraph $H$
of $G$, if $G$ is edge-colored and $H$ contains no two edges of
the same color, then $H$ is called a rainbow subgraph of $G$ and
we say that $G$ contains rainbow $H$. Let $f(G,H)$ denote the
maximum number of colors in an edge-coloring of $G$ with no
rainbow $H$. Define $rb(G,H)$ the minimum number of colors such
that any edge-coloring of $G$ with at least $rb(G,H)=f(G,H)+1$
colors contains a rainbow subgraph $H$. $rb(G,H)$ is called the
{\it rainbow number} of $H$ with respect to $G$.

When $G=K_n$, $f(G,H)$ is called the {\it anti-Ramsey number} of
$H$. Anti-Ramsey numbers were introduced by Erd\H{o}s, Simonovits
and S\'{o}s in the 1970s. Let $P_k$ and $C_k$ denote the path and
the cycle with $k$ edges, respectively. Simonovits and S\'{o}s
\cite{Simon} determined $f(K_n,P_k)$ for large enough $n$.
Erd\H{o}s et al. \cite{Erd} conjectured that for every fixed
$k\geq 3$,
$f(K_n,C_k)=n(\frac{k-2}{2}+\frac{1}{k-1})+\textit{O}(1)$, and
proved it for $k=3$ by showing that $f(K_n,C_3)=n-1$. Alon
\cite{Alon} showed that $f(K_n,C_4)=\lfloor\frac{4n}{3} \rfloor
-1$, and the conjecture is thus proved for $k=4$. Recently the
conjecture is proved for all $k\geq 3$ by Montellano-Ballesteros
and Neumann-Lara \cite{Monte}.  Axenovich, Jiang and K\"{u}ndgen
 \cite{Axeno} determined $f(K_{m,n},C_{2k})$ for all $k\geq 2$.

In 2004, Schiermeyer \cite{Schier} determined the rainbow numbers
$rb(K_n,K_k)$ for all $n\geq k\geq 4$, and the rainbow numbers
$rb(K_n,mK_2)$ for all $m\geq 2$ and $n\geq 3m+3$, where $mK_2$ is
a matching of size $m$. Li, Tu and Jin \cite{Li1} proved that
$rb(K_{m,n},pK_2)=m(p-2)+2$ for all $m\geq n\geq p\geq 3$. Chen,
Li and Tu  \cite{Li2} determined $rb(K_{n},mK_2)$.

Let $B_{n,k}$ be a $k$-regular bipartite graph with bipartition
$(X,Y)$ such that $|X|=|Y|=n$ and $k\leq n$. In this paper we give
an upper and lower bound for $rb(B_{n,k},mK_2)$, and show that for
given $k$ and $m$, if $n$ is large enough, $rb(B_{n,k},mK_2)$ can
reach the lower bound. We also determine the rainbow numbers of
matchings in paths and cycles.

\section{\large Rainbow number of matchings in regular bipartite graphs}

Denote by $mK_2$ a matching of size $m$ and $B_{n,k}$ a
$k$-regular bipartite graph with bipartition $(X,Y)$ such that
$|X|=|Y|=n$ and $k\leq n$. From a result of Li, Tu and Jin in
\cite{Li1} we know that if $n\geq 3$ and $2\leq m \leq n$, then
$rb(B_{n,n},mK_2)=n(m-2)+2$. In this section we discuss the
rainbow number of matchings in a $k$-regular bipartite graph
$B_{n,k}$.

A vertex cover of $G$ is a set $S$ of vertices such that $S$
contains at least one end-vertex of every edge of $G$. For any
$U\subset V(G)$, denote by $N_G (U)$ the neighborhood of $U$ in
$G$, abbreviate it as $N(U)$ when there is no ambiguity.

\begin{lemma}
For any bipartite graph $G$, the size of a maximum matching equals
the size of a minimum vertex cover. Let $P$ be a minimum vertex
cover of $G$, then every maximum matching of $G$ saturates $P$.
\end{lemma}

\begin{lemma}
Let $B$ be a bipartite graph with bipartition $(X,Y)$ and $|X|\geq
|Y|$, $M$ be a maximum matching of $B$. Then there exists an
$S\subseteq X$ such that $|X|-|M|=|S|-|N_B(S)|$ and $M$ saturates
$N_B(S)\cup (X-S)$, moreover, $N_B(S)\cup (X-S)$ is a minimum
vertex cover of $B$.
\end{lemma}

Let $ext(G,H)$ denote the maximum number of edges that $G$ can have
with no subgraph isomorphic to $H$.

\begin{theorem}
For any subgraph $G$ of $B_{n,k}$, if $|E(G)|>k(m-1)$, $2\leq m\leq
n$, then $mK_2\subset G$. That is
$$ext(B_{n,k},mK_2)=k(m-1).$$
\end{theorem}

{\flushleft\it Proof.}\quad By contradiction. Suppose $G$ is a
subgraph of $B_{n,k}$ with $|E(G)|>k(m-1)$ and contains no $mK_2$.
Then $G$ is bipartite and the maximum degree of the vertices in $G$
is $k$. By Lemma 2.1 $G$ has a vertex cover of size at most $m-1$,
which can cover at most $(m-1)k$ edges, contrary to $|E(G)|>k(m-1)$.
\qed

\begin{theorem}
For any $B_{n,k}$, $1\leq m\leq n$,
$$k(m-2)+2\leq rb(B_{n,k},mK_2) \leq k(m-1)+1.$$
\end{theorem}

{\flushleft\it Proof.}\quad The upper bound is obvious from Theorem
2.3. For the lower bound, let $B_{n,k}=(X,Y)$ and $Y_1\subset Y$
with $|Y_1|=m-2$, color the $k(m-2)$ edges between $Y_1$ and $X$
with $k(m-2)$ distinct colors and the remaining edges with one extra
color. It is easy to check that $k(m-2)+1$ colors are used and there
is no rainbow $mK_2$ in such a coloring. \qed

The following theorem shows that for given $k$ and $m$, if $n$ is
large enough, $rb(B_{n,k},mK_2)$ will always be equal to the lower
bound $k(m-2)+2$.

\begin{theorem}
For any given $m\geq 2$, $k\geq 3$, if $n>3(m-1)$, then
$$rb(B_{n,k},mK_2)=k(m-2)+2.$$
\end{theorem}

{\flushleft\it Proof.}\quad From Theorem 2.4 it suffices to show
that for any $m\geq 2$, $k\geq 3$, if $n>3(m-1)$, any coloring $c$
of $B_{n,k}$ with $k(m-2)+2$ colors contains a rainbow $mK_2$. By
contradiction, suppose there is no rainbow $mK_2$ in $B_{n,k}$.
Let $G$ be a subgraph of $B_{n,k}$ formed by taking one edge of
each color from $B_{n,k}$. We have $|E(G)|=k(m-2)+2$ and there is
no $mK_2$ in $G$. If there are two edge-disjoint matchings of size
$m-1$, say $M$ and $M'$ in $G$ and there exists an edge $e$ in
$B_{n,k}$ which is independent of all the edges in $M\cup M'$,
without loss of generality, say $c(e)\in C(M)$, then $M' \cup \{e
\}$ is a rainbow $mK_2$ in $B_{n,k}$. So we now focus on $G$ and
will first prove that there are two edge-disjoint matchings of
size $m-1$ in $G$.

We claim that there exists a matching of size $m-1$ saturating all
the vertices of degree $k$ in $G$. Since $|E(G)|=k(m-2)+2>k(m-2)$,
by Theorem 2.3 there is at least one $(m-1)K_2$ in $G$. Let $M$ be
a maximum matching of size $m-1$ containing maximum number of
vertices of degree $k$ in $G$, denote by $M_X$ and $M_Y$ the sets
of vertices covered by $M$ in $X$ and $Y$, respectively. If $M$
saturates all the vertices of degree $k$, we are done. Otherwise
let $v$ be a vertex of degree $k$ that is not covered by $M$,
without loss of generality, let $v\in X$. From the maximality of
$M$, every vertex in $N(v)\subset Y$ is saturated by $M$. Denote
by $MN(v)$ the set of vertices corresponding to $N(v)$ through
$M$, that is $MN(v)=\{ u | uv'\in M, v'\in N(v) \}$. If there is a
vertex $v_1\in MN(v)$ with $d(v_1)<k$ corresponding to $u_1\in
N(v)$ in $M$, let $M_1=M\cup \{ u_1v \} \backslash \{ u_1v_1 \}$
and we get a matching of size $m-1$ which has more vertices of
degree $k$ than $M$, contrary to the choice of $M$. So every
vertex in $MN(v)$ is of degree $k$. Since the sum of the degrees
of vertices in $N(v)$ is at most $k|N(v)|$ and the sum of the
degrees of vertices in $MN(v)\cup \{ v \}$ is $(k+1)|N(v)|$, and
there is no augmenting path in $G$ with respect to $M$, there
exists an edge $v_2u_3$ with $v_2\in MN(v)$ and $u_3\in M_X
\backslash N(v)$. Let $u_2$ and $v_3$ be the corresponding
vertices of $v_2$ and $u_3$ in $M$, respectively. If $d(v_3)<k$,
then $M_2=M\cup \{ u_2v,u_3v_2 \} \backslash \{ u_2v_2,u_3v_3 \}$
is a matching with more vertices of degree $k$ than $M$, a
contradiction. So $d(v_3)=k$. Now the sum of the degrees of
vertices in $N(v)\cup \{u_3\}$ is at most $(k+1)|N(v)|$ and the
sum of the degrees of vertices in $MN(v)\cup \{ v,v_3 \}$ is
$(k+2)|N(v)|$. There is an edge $v_4u_5$ with $v_4\in MN(v)\cup \{
v_3 \}$ and $u_5\in M_X \backslash N(v)$. Continue this procedure
recursively. Since there are finite vertices in $M_X$ and there is
at least one vertex of degree less than $k$ in $M_X$, we can stop
at a vertex $w\in M_X$ with $d(w)<k$ and get a matching $M'$ with
vertex set $V(M')=V(M)\backslash \{w\} \cup \{v\}$, contrary to
the choice of $M$.

Let $M$ be a matching of size $m-1$ saturating all the vertices of
degree $k$ in $G$ and $G'=G\backslash M$. Then the maximum degree
of the vertices in $G'$ is $k-1$ and
$|E(G')|=k(m-2)+2-m+1=km-2k-m+3>(k-1)(m-2)$, and so the size of
the minimum vertex cover is at least $m-1$. By Lemma 2.1 there is
a matching $M'$ of size $m-1$ in $G'$. Now $M$ and $M'$ are two
edge-disjoint matchings of size $m-1$ in $G$.

Since $M$ and $M'$ are both maximum matchings in $G$, by Lemma 2.1
the edges in $M\cup M'$ are incident to at most $3(m-1)$ vertices,
which can be incident to at most $3k(m-1)$ edges. If $n>3(m-1)$,
then $|E(B_{n,k})|>3k(m-1)$ and there is at least one edge in
$B_{n,k}$ that is independent of $E(M)\cup E(G')$, which completes
the proof. \qed

\section{\large Rainbow numbers of matchings in paths and cycles}

In this section we suppose $n\geq 3$. Let $P_n$ be the path with
$n$ edges with $V(P_n)=\{ x_0,x_1,\cdots,x_n\}$ and $E(P_n)=\{e_i|
e_i=x_{i-1}x_i,1\leq i \leq n\}$, and let $C_n$ be the cycle with
$n$ edges.

\begin{theorem}
For any $1\leq m\leq \lceil \frac{n}{2}\rceil $,

$$2m-2\leq rb(P_n,mK_2) \leq 2m-1.$$
\end{theorem}

{\flushleft\it Proof.}\quad For the upper bound, let $c$ be any
coloring of $P_n$ with $2m-1$ colors, and $G$ be the spanning
subgraph formed by taking one edge of each color from $P_n$. Then
$G$ is a bipartite graph, and so the size of its maximum matchings
equals the size of its minimum vertex covers. Since one vertex can
cover at most two edges in $G$, the size of a minimum vertex cover
of $G$ is at least $m$, and so there is a matching of size $m$ in
$G$ and hence there is a rainbow $mK_2$ in $P_n$.

To obtain the lower bound we need to show that there is a coloring
$c$ of $P_n$ with $2m-3$ colors without rainbow $mK_2$. Let
$c(e_i)=i$ for $i=1,\cdots,2m-4$ and color all the other edges
with $2m-3$. It is easy to see that there is no rainbow $mK_2 $ in
such a coloring. \qed

The following theorem gives a relationship between $rb(G,mK_2)$
and $rb(H,mK_2)$, in which $H$ is obtained from $G$ by identifying
two vertices of $G$ without any common neighbor.

\begin{theorem}
Let $G$ be a graph, $x',x''\in V(G)$ with $N(x')\cap
N(x'')=\emptyset$. Identify $x'$ and $x''$ into one vertex $x$ and
let the resultant graph be $H$, that is $V(H)=V(G)\cup \{ x\}
\backslash \{x',x''\}$ and $E(H)=\{ uv| uv\in E(G)\ \text{and} \
\{ u,v\} \cap \{x',x''\}= \emptyset \} \cup \{ xu|x'u\in E(G)\
\text{or} \ x''u\in E(G)\}$. Then $rb(G,mK_2)\leq rb(H,mK_2)$.
\end{theorem}

{\flushleft\it Proof.}\quad Let $rb(H,mK_2)=p$ and $c$ be any
coloring of $G$ with $p$ colors. For each edge in $G$, color the
corresponding edge in $H$ with the same color. Then there is a
rainbow $mK_2$ in $H$. Since the corresponding edge set in $G$ of
an independent edge set in $H$ is still independent, we have a
rainbow $mK_2$ in $G$, and so $rb(G,mK_2)\leq p$. \qed

Notice that $C_n$ can be obtained from $P_n$ by identifying the
two ends of $P_n$. Thus from above theorem we have

\begin{corollary}
$rb(P_n,mK_2)\leq rb(C_n,mK_2)$.
\end{corollary}

In Theorem 3.1, if we replace $P_n$ by $C_n$ and $m\leq \lceil
\frac{n}{2}\rceil $ by $m\leq \lfloor \frac{n}{2}\rfloor $, then
from Corollary 3.3 we get the following theorem.

\begin{theorem}
For any $1\leq m\leq \lfloor \frac{n}{2}\rfloor $,
$$2m-2\leq rb(C_n,mK_2) \leq 2m-1.$$
\end{theorem}

\begin{theorem}
For any $2 \leq m\leq \lceil \frac{n}{2}\rceil $,
$$rb(P_n,mK_2)=\{ \begin{array}{cc}
                 2m-1, & n\leq 3m-3; \\
                 2m-2, & n>3m-3.
               \end{array}
 $$
\end{theorem}

{\flushleft\it Proof.}\quad For $n\leq 3m-3$, since $2m-2\leq
rb(P_n,mK_2) \leq 2m-1$, we can construct a coloring of $P_n$ with
$2m-2$ colors that contains no rainbow $mK_2$. In fact, let
$p=n-(2m-2)$, and for $1\leq i\leq p$ let
$c(e_{3i-2})=c(e_{3i})=2i$ and $c(e_{3i-1})=2i-1$, and for $1\leq
j\leq n-3p$ let $c(e_{3p+j})=2p+j$. It is easy to check that for
such a coloring, in any rainbow matching of $P_n$ only one color
of $2i-1$ and $2i$ ($1\leq i\leq m-1$) may appear, and so there is
no rainbow $mK_2$ in $P_n$.

For $n>3m-3$, let $c$ be any coloring of $P_n$ with $2m-2$ colors.
We will prove that there is a rainbow $mK_2$ in $P_n$. By
contradiction, suppose there is no rainbow $mK_2$ in $P_n$. Let
$G$ be the spanning subgraph of $P_n$ formed by taking one edge of
each color in $P_n$, $E(G)=\{e_{i_1},e_{i_2},\cdots,
e_{i_{2m-2}}\}$, $ 1\leq i_1 < i_2 < \cdots < i_{2m-2} \leq n $
with $c(e_{i_j})=j,\ 1\leq j \leq 2m-2 $. There is no $mK_2$ in
$G$. Notice that $G$ is bipartite, and so the size of maximum
matchings equals the size of minimum vertex covers. Since one
vertex of $G$ can cover at most two edges, there is a vertex cover
of size $m-1$ in $G$, and so
$e_{i_{2l-1}}$ is adjacent to $e_{i_{2l}},\ 1\leq l \leq m-1 $. \\

\noindent{\it Claim 1.} Every edge $e$ in $P_n\backslash E(G)$ is
adjacent to an edge in $E(G)$. Otherwise suppose there is an edge
$e\in E(P_n) \backslash E(G)$ independent of $E(G)$. Notice that
$M_{1}=\{e_{i_1},e_{i_3}, \cdots, e_{i_{2m-3}} \}$ and
$M_{2}=\{e_{i_2},e_{i_4},\cdots, e_{i_{2m-2}} \}$ are two disjoint
matchings of size $m-1$ in $G$. Let $c(e)=c(e_{i_l})$, and without
loss of generality, let $e_{i_l}\in M_1$. Then $M_2 \cup \{ e \}$
is a rainbow $mK_2$ in $P_n$, a contradiction.\\

\noindent{\it Claim 2.}  There is no subgraph isomorphic to $P_3$ in
$P_n\backslash E(G)$. Otherwise the middle edge of $P_3$ is
independent of $E(G)$, which is contrary to Claim 1.\\

From Claims 1 and 2 we know that every nontrivial component of
$P_n\backslash E(G)$ is a single edge $P_1$ or a $P_2$. We
consider three cases and each leads to a contradiction.\\

\noindent{\it Case 1.} All the nontrivial components of
$P_n\backslash E(G)$ are single edges. From Claim 1 and $n>3m-3$,
we can deduce that $n=3m-2$ and
$E(G)=\{e_{2},e_{3},e_{5},e_{6},e_{8},e_{9}$, $\cdots,
e_{3m-4},e_{3m-3} \}$ with $c(e_{3i-1})=2i-1,\ c(e_{3i})=2i,\
1\leq i \leq  m-1$. Now $M_{1}^{1} =\{e_{3i} | 1\leq i \leq m-1
\}$ and $M_{2}^{1} =\{e_3 \} \cup \{e_{3i-1} | 2\leq i \leq m-1
\}$ have only $e_3$ in common and both are independent of $e_1$.
To avoid the existence of a rainbow $mK_2$ in $P_n$, we have
$c(e_1)=c(e_3)=2$. Similarly, $M_{1}^{2} =\{e_1,e_6 \}
\cup\{e_{3i-1} | 3\leq i \leq m-1 \}$ and $M_{2}^{2} =\{e_2,e_6 \}
\cup \{e_{3i} | 3\leq i \leq m-1 \}$ have only $e_6$ in common and
both are independent of $e_4$, and $c(e_4)=c(e_6)=4$. By the same
method, we know that $c(e_{3i-2})=c(e_{3i})=2i, 1\leq i \leq m-1$.
Then, $M_{1}^{m} =\{e_{3i-2} | 1\leq i \leq m-1 \}$ and $M_{2}^{m}
=\{e_{3i-1} | 1\leq i \leq m-1 \}$ are disjoint and both are
independent of $e_{3m-2}$. Whatever color $e_{3m-2}$ receives, we
will get a rainbow $mK_2$ in $P_n$, a contradiction.\\

Now at least one component of $P_n\backslash E(G)$ is isomorphic
to $P_2$.\\

\noindent{\it Case 2.} At least one of the end edges of $P_n$ is
in $P_n\backslash E(G)$. Without loss of generality, let
$E(G)=\{e_{2},e_{3},e_{6},e_{7},e_{9},e_{10}$, $\cdots,
e_{3m-3},e_{3m-2} \}$ with $c(e_2)=1$, $c(e_3)=2$,
$c(e_{3i})=2i-1,\ c(e_{3i+1})=2i,\ 2\leq i \leq  m-1$. Since
$M_{1}' =\{e_{3i} | 1\leq i \leq m-1 \}$ and $M_{2}' =\{e_3 \}
\cup \{e_{3i+1} | 2\leq i \leq m-1 \}$ have only $e_3$ in common
and both are independent of $e_1$, $c(e_1)=c(e_3)=2$. Now $M_{1}''
=\{e_1\} \cup\{e_{3i} | 2\leq i \leq m-1 \}$ and $M_{2}'' =\{e_2\}
\cup \{e_{3i+1} | 2\leq i \leq m-1 \}$ are disjoint and both are
independent of $e_4$. Whatever color $e_4$ receives, we will
get a rainbow $mK_2$ in $P_n$.\\

\noindent{\it Case 3.} Since none of the end edges of $P_n$ is in
$P_n\backslash E(G)$, there are at least two components in
$P_n\backslash E(G)$ isomorphic to $P_2$. Without loss of
generality, let
$E(G)=\{e_{1},e_{2},e_{5},e_{6},e_{9},e_{10},e_{12},e_{13},
\cdots, e_{3m-3},e_{3m-2} \}$ with $c(e_2)=1$, $c(e_2)=2$,
$c(e_5)=3$, $c(e_6)=4$, $c(e_{3i})=2i-1,\ c(e_{3i+1})=2i,\ 3\leq i
\leq  m-1$. Since $M_{1}' =\{e_1,e_3 \} \cup \{e_{3i} | 3\leq i
\leq m-1 \}$ and $M_{2}' =\{e_1,e_4 \} \cup \{e_{3i+1} | 2\leq i
\leq m-1 \}$ have only $e_1$ in common and both are independent of
$e_3$, we have $c(e_1)=c(e_3)=1$. $M_{1}'' =\{e_1,e_6 \} \cup
\{e_{3i} | 3\leq i \leq m-1 \}$ and $M_{2}'' =\{e_2,e_6 \} \cup
\{e_{3i+1} | 2\leq i \leq m-1 \}$ have only $e_6$ in common and
both are independent of $e_4$, we have $c(e_4)=c(e_6)=4$. Now
$M_{1} =\{e_1,e_4 \} \cup \{e_{3i} | 3\leq i \leq m-1 \}$ and
$M_{2} =\{e_2,e_5 \} \cup \{e_{3i+1} | 2\leq i \leq m-1 \}$ are
disjoint and both are independent of $e_7$. Whatever color $e_7$
receives, we will get a rainbow $mK_2$ in $P_n$.\qed

From Corollary 3.3 and Theorem 3.5, we have $rb(C_n,mK_2)=2m-1,\
n\leq 3m-3$. For $n>3m-3$, by a similar proof in Theorem 3.5, we
have $rb(C_n,mK_2)=2m-2$. Thus we have

\begin{theorem}
For any $m\leq \lfloor \frac{n}{2}\rfloor $,
$$rb(C_n,mK_2)=\{ \begin{array}{cc}
                 2m-1, & n\leq 3m-3; \\
                 2m-2, & n>3m-3.
               \end{array}
 $$
\end{theorem}

\end{document}